\documentclass[12pt,leqno]{amsart}

\usepackage{tikz}
\usepackage{tikz-cd}
\tikzset{>=stealth}
\tikzcdset{arrow style=tikz}

\usepackage{enumitem}
\setlist[enumerate,1]{label=\textnormal{(\roman*)}}
\setlist[itemize,1]{label={--}}

\usepackage{amsxtra}
\usepackage{amsmath}
\usepackage{amstext}
\usepackage{amssymb}
\usepackage{mathrsfs}

\usepackage{pdfsync}
\usepackage[obeyspaces,hyphens,spaces]{url}
\usepackage[colorlinks=true,linkcolor=violet,citecolor=blue]{hyperref}

\usepackage{inconsolata}
\usepackage{fourier}
\usepackage[T1]{fontenc}

\usepackage[ left = 3cm, right = 3cm,
             top = 3cm, bottom = 3cm ]{geometry}

\makeatletter
\newlength\@oldparskip
\let\@oldtoc=\tableofcontents
\def\tableofcontents{\setlength{\@oldparskip}{\parskip}%
  \setlength{\parskip}{0pt}%
  \@oldtoc%
  \setlength{\parskip}{\@oldparskip}}

\counterwithin{paragraph}{section}
\counterwithin{equation}{subsection}
\newcommand{\epointn}[1]{\vspace{2mm}\par{\em #1.} }

\def\section{\@startsection{section}{1}%
  \z@{.7\linespacing\@plus\linespacing}{.5\linespacing}%
  {\normalfont\scshape\centering}}

\def\tagform@#1{\maketag@@@{(\ignorespaces{\bfseries#1}\unskip\@@italiccorr)}}
\def\tagformnorm@#1{\maketag@@@{(\ignorespaces#1\unskip\@@italiccorr)}}
\renewcommand{\eqref}[1]{\textup{\tagformnorm@{\ref{#1}}}}

\makeatother

\newcommand{\psqedsymb}{$\blacksquare$}

\renewcommand\qed%
   {~\hfill{\psqedsymb}}

\newcommand\pf[1][]{\textbf{Proof}%
  \ifthenelse{\equal{#1}{}}{}{ #1}\textbf{.}\enspace}

\newlength\stdlisttopsep
\newcommand\stdlist
  {\setlength  \labelsep  {1ex}%
   \setlength  \labelwidth{2em}%
   \addtolength\labelwidth{-\labelsep}%
   \setlength  \topsep{\stdlisttopsep}%
   \setlength  \partopsep{0pt}%
   \setlength  \parsep{\parskip}%
   \setlength  \itemsep{0pt}%
   \setlength  \leftmargin{\parindent+\labelwidth}%
   \setlength  \rightmargin{0pt}%
   \setlength  \itemindent{0pt}%
   \setlength  \listparindent{\parindent}}

  {\begin{list}{{\fontfamily{cmr}\selectfont\textbullet}}%
  {\stdlist}}%
  {\end{list}}

  {\begin{list}{{\fontfamily{cmr}\selectfont\textbullet}}{\stdlist}}%
  {\end{list}}

\newcounter{enumctr}
\renewcommand\theenumctr{\roman{enumctr}}

  {\begin{list}{}%
     {\stdlist\setcounter{enumctr}{1}%
      \renewcommand\makelabel%
         {(\theenumctr\refstepcounter{enumctr})\hfil}}}%
  {\end{list}}

\newcommand{\mynewtheorem}[2]{%
  \newtheorem{#1}[subsection]{#2}%
  \newtheorem*{#1*}{#2}%
  \newtheorem{#1sub}[subsubsection]{#2}}

\mynewtheorem{thm}{Theorem}
\mynewtheorem{cor}{Corollary}
\mynewtheorem{conj}{Conjecture}
\mynewtheorem{lem}{Lemma}
\mynewtheorem{prop}{Proposition}
\mynewtheorem{clm}{Claim}

\mynewtheorem{theorem}{Theorem}
\mynewtheorem{corollary}{Corollary}
\mynewtheorem{conjecture}{Conjecture}
\mynewtheorem{example}{Example}
\mynewtheorem{lemma}{Lemma}
\mynewtheorem{proposition}{Proposition}
\mynewtheorem{claim}{Claim}

\theoremstyle{definition}
\mynewtheorem{defn}{Definition}
\mynewtheorem{notn}{Notation}
\mynewtheorem{eg}{Example}
\mynewtheorem{question}{Question}

\theoremstyle{remark}
\mynewtheorem{rem}{Remark}
\mynewtheorem{remark}{Remark}

\usepackage{mathtools}
\usepackage{stmaryrd}
\usepackage{amssymb}
\usepackage{stackrel}
\usepackage{braket}

\DeclarePairedDelimiter{\pairing}{\langle}{\rangle}

\DeclareMathOperator{\Hom}{Hom}

\DeclareMathOperator{\Aut}{Aut}

\DeclareMathOperator{\cyc}{cyc}

\newcommand{\FF}{\mathbf{F}}

\newcommand{\RR}{\mathbf{R}}

\newcommand{\ZZ}{\mathbf{Z}}

\newcommand{\sJ}{\mathscr{J}}

\newcommand{\rmC}{\mathrm{C}}

\newcommand{\rmH}{\mathrm{H}}

\newcommand{\wc}{{\mkern 2mu\cdot\mkern 2mu}}

\makeatletter
\newcommand*{\coloneqq}{\mathrel{\rlap{%
           \raisebox{0.3ex}{$\m@th\cdot$}}%
           \raisebox{-0.3ex}{$\m@th\cdot$}}%
           =}
\newcommand*{\eqqcolon}{=%
           \mathrel{\rlap{%
           \raisebox{0.3ex}{$\m@th\cdot$}}%
           \raisebox{-0.3ex}{$\m@th\cdot$}}}
\makeatother

\newcommand{\multichoose}[2]{\ensuremath{\left(\kern-.3em\left(\genfrac{}{}{0pt}{}{#1}{#2}\right)\kern-.3em\right)}}

\newcommand{\rmsup}[2]{#1^{\mathrm{#2}}}

\newcommand{\ssth}{\textsuperscript{th}}

\usepackage{graphicx}
\usepackage{caption}
\usepackage{subcaption}
\usetikzlibrary{decorations.markings, calc, positioning}

\DeclareMathOperator{\Sh}{Sh}

\DeclareMathOperator{\gr}{gr}

\newcommand{\tint}{{\textstyle\int}}

\title[Combinatorial Iterated Integrals and the Harmonic Volume]{Combinatorial Iterated Integrals and the Harmonic Volume of Graphs}
\author{Raymond Cheng}
\address{Department of Mathematics \\
  Columbia University \\
  2990 Broadway \\
  New York, NY, 10027, USA
}
\email{rcheng@math.columbia.edu}
\urladdr{http://math.columbia.edu/~rcheng/}
\author{Eric Katz}
\address{Department of Mathematics \\
  The Ohio State University \\
  231 West 18\textsuperscript{th} Avenue \\
  Columbus, OH, 43210, USA}
\email{katz.60@osu.edu}
\urladdr{https://people.math.osu.edu/katz.60/}
\date{\today}


\begin{document}
\begin{abstract}
  Let $\Gamma$ be a connected bridgeless metric graph, and fix a point $v$ of $\Gamma$.
  We define combinatorial iterated integrals on $\Gamma$ along closed paths at $v$, a unipotent generalization of the usual cycle pairing
  and the combinatorial analogue of Chen's iterated integrals on Riemann
  surfaces. These descend to a bilinear pairing between the group algebra of the fundamental group of $\Gamma$ at $v$ and
  the tensor algebra on the first homology of $\Gamma$,
  $\int\colon \ZZ\pi_1(\Gamma,v) \times T\rmH_1(\Gamma,\RR) \to \RR$.
  We show that this pairing on the two-step unipotent quotient of the group algebra allows one to
  recover the base-point $v$ up to well-understood finite ambiguity. We encode the data of this structure
  as the combinatorial harmonic volume which is valued in the tropical intermediate Jacobian. We also give a potential-theoretic
  characterization for hyperelliptiicity for graphs.
\end{abstract}
\maketitle
\setcounter{tocdepth}{1}


\setcounter{section}{0}
\section{Introduction}\label{S:intro}

A metric graph $\Gamma$ is equipped with a cycle pairing $\pairing{\wc,\wc}$ on its homology $\rmH_1(\Gamma,\RR)$ taking a pair of cycles to the signed length of their  intersection. In tropical geometry, metric graphs are the analogues of curves. Under this analogy, the cycle pairing corresponds to the integrals of holomorphic $1$-forms along closed paths. Many authors~\cite{CV10,MZ-Jacobian} have considered a setting in which the cycle pairing is packaged
as a tropical Jacobian. Indeed, by the tropical Torelli theorem~\cite{Ger82,Art06,DSV09,CV10,SW10}, the tropical Jacobian determines a bridgeless finite connected metric graph up to an equivalence relation called $2$-isomorphism.

A natural generalization of period integrals is  Chen's theory of iterated integrals~\cite{Che77}
where one integrates a number of $1$-forms along a path. This theory interpolates between homology and the fundamental
group of an algebraic curve. In fact, it encodes information about the unipotent fundamental group which can be
understood in terms of a particular quotient of the group algebra of the fundamental group. In this paper,
we introduce \emph{combinatorial iterated integrals} where we integrate
a number of homology classes $\omega_1,\dots,\omega_{\ell}\in \rmH_1(\Gamma,\RR)$ along a path $\gamma$ in $\Gamma$ to obtain $\int_{\gamma} \omega_1 \dots\omega_{\ell}$.
Combinatorial iterated integrals can naturally be interpreted as unipotent invariants of graphs. Let $\ZZ\pi_1(\Gamma,v)$ be the group algebra of the fundamental group of $\Gamma$ at some base-point $v$. If $J$ is the augmentation ideal of $\ZZ\pi_1(\Gamma,v)$, combinatorial iterated integration of length at most $\ell$ descends to a bilinear pairing
\[ \int \colon \ZZ\pi_1(\Gamma,v)/J^{\ell + 1} \times T_\ell \rmH_1(\Gamma,\RR) \to \RR \]
where $T_\ell\rmH_1(\Gamma,\RR)$ is the truncated tensor algebra
$T_\ell\rmH_1(\Gamma,\ZZ) \coloneqq \bigoplus_{k = 0}^\ell \rmH_1(\Gamma,\ZZ)^{\otimes k}$.
Because $\rmH_1(\Gamma,\RR)$ can be recovered from $\ZZ\pi_1(\Gamma,v)/J^{\ell + 1}$, combinatorial iterated integrals can be viewed as a structure on $\ZZ\pi_1(\Gamma,v)/J^{\ell + 1}$. We call this structure an {\em integration algebra}. As is the case in tropical geometry, this structure is a combinatorial shadow of its classical analogue. Indeed, the
pair $(\ZZ\pi_1(\Gamma,v)/J^3,\tint)$ can be interpreted as lowest weight component of the mixed Hodge structure on the unipotent fundamental group of a degenerating family of algebraic curves.

It is natural to ask if the graph can be recovered from $(\ZZ\pi_1(\Gamma,v)/J^3,\tint)$. Indeed, its classical analogue is the
mixed Hodge structure on truncations of the fundamental group algebra of complex algebraic varieties introduced by Hain~\cite{Hai-dR-I,Hai-dR-II,Hai85}, generalizing work of Morgan~\cite{Mor78}. Hain~\cite{Hai85} and Pulte~\cite{Pul88}, drawing on work by Carlson~\cite{Car-Ext,Car-Bowdoin},
use the mixed Hodge structure on $\ZZ\pi_1(X,x)/J^3$ to prove a Torelli Theorem for pointed complex algebraic curves $(X,x)$. We believe $(\ZZ\pi_1(\Gamma,v)/J^3,\tint)$ to be a complete invariant of connected bridgeless pointed metric graphs:
\begin{conjecture}[Unipotent Torelli Conjecture] \label{conj:unipotenttorelli}
Let $(\Gamma,v)$ be a connected bridgeless pointed metric graph.
  Then the pair $(\ZZ\pi_1(\Gamma,v)/J^3,\tint)$ completely determines $(\Gamma,v)$.
\end{conjecture}
Unfortunately, we have been unable to prove this theorem. The difficulty comes from a limitation of the tropical Torelli theorem: it only recovers a graph up to an equivalence relation called $2$-isomorphism which is generated by two moves, vertex-cleaving and Whitney twists. However, if the graph is known, one can recover the base-point of the graph from the integration algebra
up to well-understood ambiguity:
\begin{theorem} Let $\Gamma$ be a bridgeless metric graph with $g(\Gamma)\geq 2$. Let $v$ be a  point of the underlying topological space $|\Gamma|$. Then the isomorphism type of the integration algebra $\big(\ZZ\pi_1(\Gamma,v)/J^3,\tint\big)$ determines $v$ up to at most
$\left|\Aut_{\cyc}\big(\ZZ\pi_1(\Gamma,v)/J^3,\tint\big)\right|$ choices.
\end{theorem}
Here, $\Aut_{\cyc}\big(\ZZ\pi_1(\Gamma,v)/J^3,\tint\big)$ is a finite group of order at most $2^{\kappa}$ where $\kappa$ is the number of $2$-connected components of $\Gamma$.
In the case where the graph is $3$-connected, there is no ambiguity coming from $2$-isomorphism and we can recover the
graph and base-point (up to two choices) from the integration algbera:
\begin{corollary} Let $(\Gamma,v)$ be a pointed metric graphs such that $\Gamma$ is $3$-connected. Then the integration algebra $\big(\ZZ\pi_1(\Gamma,v)/J^3,\tint\big)$ determines $\Gamma$ up to tropical equivalence and determines $v$ up to two possibilities.
\end{corollary}
Here, tropical equivalence is the equivalence relation generated by subdividing edges and its inverse.

The hope of the unipotent Torelli conjecture is that pointed graphs are encoded by (a truncation of) their fundamental groups
and, therefore, are \emph{anabelian} combinatorial objects  in analogy with Grothendieck's anabelian
program in algebraic geometry~\cite{Esquisse}.

The pair $(\ZZ\pi_1(\Gamma,v)/J^3,\tint)$ can be considered as a sort of extension of $(\rmH_1(\Gamma,\ZZ),\pairing{\wc,\wc})$. The extension data can be encoded as an element of a real torus,$\sJ_2(\Gamma)$ analogous to an intermediate Jacobian in algebraic geometry. This tropical intermediate Jacobian is the recipient of an invariant $\nu_\Gamma$ analogous to the harmonic volume \cite{HarVol}. Given a connected bridgeless graph $\Gamma$, let $W(\Gamma)$ be the set of rigged
graphs $2$-isomorphic to $\Gamma$, that is triples $(\Gamma',v',\phi')$ where $(\Gamma',v')$ is a pointed metric, bridgeless graph, and $\phi'\colon\rmH_1(\Gamma,\ZZ)\to \rmH_1(\Gamma',\ZZ)$ is an isometry. Then the harmonic volume is a map
 \[  \nu_\Gamma \colon W(\Gamma)  \to \sJ_2(\Gamma). \]
Conjecturally, $\nu_\Gamma$ is an injection. However, we only have the weaker result:
\begin{theorem} Let $(\Gamma_1,v_1,\phi_1),(\Gamma_2,v_2,\phi_2)\in W(\Gamma)$. There is equality of harmonic volumes, $\nu_\Gamma(\Gamma_1,v_1,\phi_1)=\nu_\Gamma(\Gamma_2,v_2,\phi_2)$ if and only if  $\phi_2\circ\phi_1^{-1}\colon \rmH_1(\Gamma_1,\ZZ)\to \rmH_1(\Gamma_2,\ZZ)$ lifts to an isomorphism
\[\phi\colon \big(\ZZ\pi_1(\Gamma_1,v_1)/J_1^3,\tint\big)\to \big(\ZZ\pi_1(\Gamma_2,v_2)/J_2^3,\tint\big).\]
\end{theorem}

Our approach to the unipotent Torelli conjecture is blocked by automorphisms of $(\ZZ\pi_1(\Gamma,v)/J^3,\tint)$ that
preserve the cycles in $\rmH_1(\Gamma,\ZZ)$. These automorphisms prevent one from gluing isomorphisms in
an inductive argument. We conjecture that the only automorphisms of the integration algebra come from legitimate
automorphisms of a graph:
\begin{conjecture} \label{conj:hyperelliptic} For a $2$-connected graph, the group $\Aut_{\cyc}\big(\ZZ\pi_1(\Gamma,v)/J^3,\tint\big)$ is nontrivial if and only if $\Gamma$ is a hyperelliptic graph and $v$ is a fixed point of a hyperelliptic involution.
\end{conjecture}
Because hyperelliptic graphs are so central to our story, we prove Theorem~\ref{t:hyperellipticity}, a potential-theoretic criterion for hyperellipticity that was part of an unsuccessful attempt to prove Conjecture~\ref{conj:hyperelliptic}.

There are a number of questions that we would like to consider in the future. One should clarify the connection between our
integration algebras and the tropical Ceresa classes of Corey, Ellenberg, and Li \cite{CEL20}. Also, the discrete geometric
picture is lacking. The Torelli theorem for graphs was proved using Delaunay cells. Is there some unipotent analogue of
Delaunay cells that will allow us to reconstruct a graph?

We expect our combinatorial iterated integrals to have applications in number theory and Hodge theory. Indeed,
they have already been used by Betts and Dogra \cite{BettsDogra} in their study of the \'{e}tale fundamental
groupoid.  In work~\cite{KL:WIP} in preparation between the second named author and Daniel
Litt, it is shown that combinatorial iterated integrals mediate between the Berkovich~\cite{Ber-integrals} and Vologodsky~\cite{Vol-integrals}
notions of $p$-adic integration on curves.
This is related to work of Besser and Zerbes~\cite{BZ}, and has
applications to the non-abelian Chabauty methods of
Kim~\cite{Kim-Chabauty, Kim-survey}.
See~\cite[\S6.5]{KRZB-Utah} for a related discussion.
Also, combinatorial iterated integrals arise in asymptotics of variations of Hodge structures associated with truncated
fundamental group algebras on a semistable family of curves, analogous to the
cohomological case as explained in~\cite[Proposition 13.3]{Griff-Bul}
and~\cite[Theorem 6.6]{Sch-VHS}.

\epointn{Outline}\label{SS:intro-outline}
In \S\ref{S:grp-alg}, we review basic facts about group algebras.

In \S\ref{S:pairing}, we construct combinatorial iterated integrals and establish
their basic properties.
In particular, we formulate and prove Theorem~\ref{SS:pairing-duality}, a duality
result for combinatorial iterated integrals.

In \S\ref{S:automorphisms}, we discuss cyclic automorphism groups of $\rmH_1(\Gamma,\ZZ)$ and of integration
algebras. These are important for understanding the base-point ambiguity in the main theorems of the following sections.

In \S\ref{S:unipotent}, we prove our weaker versions of the unipotent Torelli conjecture allowing one to recover base-points
of a graph (up to ambiguity) from integration algebras.

In \S\ref{S:comp}, we define the harmonic volume invariant.

In \S\ref{S:hyperelliptic}, we review hyperelliptic graphs and prove our potential-theoretic criterion.

\epointn{Acknowledgements}%
The authors would like to thank Omid Amini, Matthew Baker, Dustin Cartwright,
Jordan Ellenberg, Joshua Greene, Daniel Litt, Sam Payne, Joseph Rabinoff, and Farbod Shokrieh for comments.


\section{Group Algebras}\label{S:grp-alg}
We review group algebras. Let $G$ be a group and $R$ be a commutative ring.
For an $R$-module $M$, let $TM\coloneqq \bigoplus_{k \geq 0} M^{\otimes k}$ be the tensor algebra
over $R$ of $M$. For $\ell\in\ZZ$, let $T_{\ell}M\coloneqq \bigoplus_{k \geq 0}^{\ell} M^{\otimes k}$ be its truncation.
The \emph{group algebra} $RG$ of $G$ with coefficients in $R$ is the
associative unital $R$-algebra with underlying $R$-module
\[ RG = \bigoplus\nolimits_{g \in G} R e_g \]
and multiplication determined on basis elements by
$r_1e_{g_1} \cdot r_2e_{g_2} \coloneqq r_1r_2e_{g_1g_2}$
for all $r_1,r_2 \in R$ and $g_1,g_2 \in G$.
The unit of $RG$ is the element $e_{1}$, where $1 \in G$ the identity of
$G$.
We write $g$ for $e_g$ and $1$ for the unit of $RG$,
$G$, and the ring $R$.

The kernel $J \coloneqq \ker(\epsilon \colon RG \to R)$ of the augmentation
homomorphism is called the \emph{augmentation ideal}.
An element of $J$ is of the form $\sum_{g \in G} a_gg$ with
$\sum_{g \in G} a_g = 0$.
So
\[
  \sum_{g \in G} a_g g
    = \sum_{g \in G \setminus \{1\}} a_g(g - 1) + \sum_{g \in G} a_g
    = \sum_{g \in G \setminus \{1\}} a_g(g - 1),
\]
showing that $J$ is a free $R$-module with basis $\set{(g - 1) | g \in G \setminus \{1\}}$.
Since both $R$ and $J$ are free $R$-modules, the  exact sequence
$0 \to J \to RG \xrightarrow{\epsilon} R \to 0$
is split and we have $RG \cong R \oplus J$ as $R$-modules.
There is an isomorphism of $R$-modules
\[\rmsup{G}{ab}  \otimes_\ZZ R \to J/J^2 \colon g \mapsto (g - 1).\]


The algebra $RG$ is has a descending filtration
$RG = J^0 \supseteq J^1 \supseteq J^2 \supseteq \cdots$
by powers of the augmentation ideal.
The associated graded algebra with respect to this filtration is
\[ \gr_J(RG) \coloneqq \bigoplus\nolimits_{i = 0}^\infty J^i/J^{i+1}. \]
After identifying $J/J^2$ with $\rmsup{G}{ab}\otimes_\ZZ R$,
there is a natural map
\[( \rmsup{G}{ab} \otimes_\ZZ R)^{\otimes i} \to J^i/J^{i+1}\colon (g_1 - 1) \otimes \cdots \otimes (g_i - 1) \mapsto (g_1 - 1) \cdots (g_i - 1) \]
for each $i \geq 0$ and where $g_1,\ldots,g_i \in \rmsup{G}{ab}$.
If $G$ is a free group on a set of generators $E$,
this map is an isomorphism and we see that the associated graded algebra of
$RG$
\[ \gr_J(RG) \cong \bigoplus\nolimits_{i = 0}^\infty (\rmsup{G}{ab}\otimes_\ZZ R)^{\otimes i} = T(\rmsup{G}{ab}\otimes_\ZZ R) \]
is isomorphic to the tensor algebra of the free $R$-module on the set $E$.

\section{Combinatorial Iterated Integrals}\label{S:pairing}

This section introduces combinatorial iterated integrals and discusses their properties.

\subsection{Background on graphs}
Throughout this paper, a \emph{graph} is a finite connected graph, possibly with loops and
multiple edges.
For topological constructions, we think of graphs as one-dimensional
$\Delta$-complexes. Given a graph $\Gamma$, we write $V(\Gamma)$ and $E(\Gamma)$ for the set of
vertices and oriented edges of $\Gamma$, respectively. Write $|\Gamma|$ for its underlying topological space.
For an oriented edge $e$, let $\overline{e}$ denote $e$ with its orientation
reversed. Given an oriented edge $e$, write $e^+$ and $e^-$ for its head and tail in the given orientation. We will often assume our graphs to be \emph{bridgeless}, that is, there does not
exist an edge whose deletion disconnects the graph.  For a positive integer $k$, a graph is \emph{$k$-connected} or
\emph{$k$-vertex connected} if it cannot be disconnected by removing fewer than $k$ vertices.

A \emph{metric graph} $(\Gamma,\ell)$ is a graph $\Gamma$ together with a length function $\ell\colon E(\Gamma)\to\RR_{>0}$ such that $\ell(\overline{e})=\ell(e)$. For each oriented edge $e$ of $\Gamma$ we fix a homeomorphism $t_e\colon |e|\to [0,\ell(e)]$ such that
$t_{\overline{e}}=\ell(e)-t_e$. A \emph{polynomial $1$-form} $\omega$ on $\Gamma$ is a choice of $1$-form
$\omega_e = p_e(t) dt$ on each oriented edge $e\cong [0,\ell(e)]$ such that $p_e(t)$ is a real polynomial and
\[p_{\overline{e}}(t)=-p_e(\ell(e)-t).\]
We interpret a polynomial $1$-form as a $1$-form on $|\Gamma|$. The degree of a polynomial $1$-form is the maximum of the degrees of the polynomials $p_e(t)$.
A polynomial $1$-form is said to be a {\em tropical $1$-form} if
\begin{enumerate}
\item each $p_e(t)$ is constant; and
\item if $v\in V(\Gamma)$ and $e_1,\dots,e_k$ are the edges adjacent to $v$, directed away from $v$, $\sum_i p_{e_i}(t)=0$.
\end{enumerate}
The second condition is also called harmonicity. See \cite{MZ-Jacobian} for more details.

For a ring $R$, let $\rmC_0(\Gamma,R)$ be the free $R$-module on $V(\Gamma)$, and let $\rmC_1(\Gamma,R)$ be the quotient of the free $R$-module on $E(\Gamma)$ by the relation $\overline{e}=-e$.
Let $\partial \colon \rmC_1(\Gamma,R) \to \rmC_0(\Gamma,R)$ be the simplicial
boundary map defined by $e \mapsto e^+ - e^-$.
The homology groups of $\Gamma$ with coefficients in $R$ are
$\rmH_0(\Gamma,R) \coloneqq \rmC_0(\Gamma,R)/\partial \rmC_1(\Gamma,R)$ and
$\rmH_1(\Gamma,R) \coloneqq \ker(\partial \colon \rmC_1(\Gamma,R) \to \rmC_0(\Gamma,R))$.
The elements of $\rmH_1(\Gamma,R)$ are referred to as \emph{cycles} of $\Gamma$. A cycle $C$ is \emph{simple} if the coefficient of each edge in $C$ is $1$, $0$, or $-1$. The set of edges for which the coefficient in $C$ is nonzero is called
the \emph{support} of $C$. A simple cycle is \emph{primitive} if its support is minimal among the cycles.
To a cycle $C\in \rmH_1(\Gamma,\RR)$, we may attach a tropical $1$-form $\omega_C$ as follows:
write $C=\sum a_e e$ for $a_e\in\RR$ and set $p_e(t)=a_e$. The map $C\mapsto \omega_C$ is an isomorphism
from $\rmH_1(\Gamma,\RR)$ to tropical $1$-forms. Henceforth we shall identify elements of $\rmH_1(\Gamma,\RR)$ with
tropical $1$-forms.
For a closed path $\gamma\colon[0,1]\to |\Gamma|$, we write $[\gamma]\in \rmH_1(\Gamma,\ZZ)$ for the underlying cycle.

We define an inner product on $\rmC_1(\Gamma,\RR)$ by
\[
  \pairing{\wc,\wc} \colon \rmC_1(\Gamma,\RR) \times \rmC_1(\Gamma,\RR) \to \RR \qquad \pairing{x,y} \coloneqq \begin{dcases*} \ell(x) & if $x = y$, \\
  -\ell(x) & if $x=\overline{y}$\\
   0 & if $x \neq y,\overline{y}$. \end{dcases*}
\]
Since $\rmH_1(\Gamma,\RR) \subseteq \rmC_1(\Gamma,\RR)$, the first homology group inherits
a bilinear form, also denoted by $\pairing{\wc,\wc}$. It takes a pair of cycles to the signed length of their intersection,
counted with multiplicity. We call this the \emph{cycle pairing}.
Because it is the restriction of the standard
Euclidean pairing, it is positive-definite and thus, nondegenerate.


\subsection{Combinatorial Iterated Integrals}\label{SS:pairing}
In this subection, we introduce combinatorial iterated integrals on graphs, a non-abelian extension of the  cycle pairing.

Our constructions are inspired by Chen's theory of iterated line
integrals~\cite{Che77} and their application to the construction of a mixed
Hodge structure on the fundamental group of an algebraic variety~\cite{Hai85}.
After completing this article, we became aware of~\cite{BKP09} in which
related definitions were made.

To define iterated integrals on paths, we first define antiderivatives on the universal cover $\tilde{\Gamma}$. A continuous function $F\colon |\tilde{\Gamma}|\to \RR$ is said to be piecewise polynomial if it restricts to each edge as a polynomial by the parameterization $t_e$. For such a function $F$, we define $dF$ to be the $1$-form whose restriction to
an edge $e$ is the differential $dF|_e$. Pick a base-point $v$ of $\tilde{\Gamma}$.
Let $\omega_1,\dots,\omega_k$ be polynomial $1$-forms on $\Gamma$ which we will identify with their pulbacks on $\tilde{\Gamma}$. The primitive
$F_{\omega_1\cdots\omega_k}$ is defined by induction on $k$. We define $F_{\omega_1}\colon\tilde{\Gamma}\to \RR$
to be the continuous piecewise polynomial function with $F_{\omega_1}(v)=0$ and $dF_{\omega_1}=\omega_1$. In
general, we define $F_{\omega_1\cdots\omega_k}$ to be the piecewise polynomial function characterized by
\begin{enumerate}
\item $F_{\omega_1\cdots\omega_k}(v)=0$, and
\item $dF_{\omega_1\cdots\omega_k}=F_{\omega_1\cdots\omega_{k-1}}\omega_k$.
\end{enumerate}
For a path $\gamma\colon[0,1]\to |\Gamma|$, we define
\[\int_\gamma \omega_1\cdots\omega_k=F_{\omega_1\cdots\omega_k}(\tilde{\gamma}(1))\]
where $\tilde{\gamma}\colon [0,1]\to \tilde{\Gamma}$ is a lift of $\gamma$ to $\tilde{\Gamma}$ and the primitive
is taken at the base-point $\tilde{\gamma}(0)$.

By extending linearly, we produce \emph{combinatorial iterated integrals} as a bilinear map
\[ \int \colon \ZZ\pi_1(\Gamma,v) \times T\rmH_1(\Gamma,\ZZ) \to \RR. \]
We follow the convention that the integral against the identity element $1\in \rmH_1(\Gamma,\ZZ)^{\otimes 0}\subset T\rmH_1(\Gamma,\ZZ)$ is the augmentation map:
\[\int_\gamma 1 = \epsilon(\gamma).\]
The following lemma follows from unwinding definitions:
\begin{lemma} Let $\gamma\colon[0,1]\to |\Gamma|$ be a loop, and let $\omega_C$ be the tropical $1$-form attached to $C\in \rmH_1(\Gamma,\RR)$. Then
\[\int_{\gamma} \omega_C=\pairing{[\gamma],C}.\]
\end{lemma}

\subsection{Properties of combinatorial iterated integrals}

Combinatorial iterated integrals can be expressed in terms of shuffles.
For positive integers $k$ and $\ell$, define the set $\Sh(k,\ell)$ of
\emph{$(k,\ell)$-shuffles} to be the following subset of the symmetric group
$S_{k+\ell}$ on $k + \ell$ symbols:
\[
  \Sh(k,\ell) \coloneqq
    \Set{\sigma \in S_{k + \ell} | \sigma^{-1}(1) \leq \cdots \leq \sigma^{-1}(k) \;\text{and}\; \sigma^{-1}(k+1) \leq \cdots \leq \sigma^{-1}(k + \ell)}.
\]

Combinatorial iterated integrals have properties analogous to that of classical iterated integrals. The proofs are analogous:
one rewrites the iterated integral as the integral of a $k$-form on a time-ordered simplex and uses properties of integration.
\begin{proposition}\label{SS:pairing-hopf}
Let $\alpha$ be a path in $\Gamma$, and let $\omega_1,\ldots,\omega_k$ be polynomial $1$-forms. Then
  we have the following formulae:
  \begin{itemize}
    \item (Product) For any
      $\omega_{k + 1},\ldots,\omega_{k + \ell} \in \rmC_1(\Gamma,\ZZ)$,
      \[
        \int_\alpha \omega_1 \cdots \omega_k \int_\alpha \omega_{k + 1} \cdots \omega_{k + \ell} = \sum_{\sigma \in \Sh(k,\ell)} \int_\alpha \omega_{\sigma(1)} \cdots \omega_{\sigma(k + \ell)}.
      \]
    \item (Concatenation) For any path $\beta$ with $\beta(0)=\alpha(1)$,
      \[
        \int_{\alpha\beta} \omega_1 \cdots \omega_k = \sum_{i = 0}^k \int_{\alpha} \omega_1 \cdots \omega_i \int_{\beta} \omega_{i+1} \cdots \omega_k.
      \]
    \item (Antipode)
      \[
        \int_{\alpha^{-1}} \omega_1 \cdots \omega_k = (-1)^k \int_\alpha \omega_k \cdots \omega_1.
      \]
  \end{itemize}
\end{proposition}

\begin{proof}
  See, for example,~\cite[2.9,2.11,2.12]{Hai85}.
\end{proof}

For any positive integer $n$, write $[n] \coloneqq \{1,\ldots,n\}$.
For $k$ and $r$ positive integers, let
\[ \Delta(k,r) \coloneqq \Set{f \colon [k] \to [r] | f(1) \leq \cdots \leq f(k)} \]
be the set of all weakly increasing functions from $[k]$ to $[r]$.
Equivalently, setting $n_i \coloneqq \#f^{-1}(i)$ for $i = 1,\ldots,r$,
an element $f \in \Delta(k,r)$ may be represented as the sequence $(n_1,\ldots,n_r)$.
Note $n_1 + \cdots + n_r = k$.

The next two formulas are obtained by iterating those proven above.

\begin{theorem} We have the following formulas:
\begin{itemize}
\item (Symmetrization Formula)\label{SS:pairing-sym-product}
\emph{Let $\gamma$ be a path in $\Gamma$, and let $\omega_1,\ldots,\omega_k$ be polynomial $1$-forms.
  Then
  \[ \sum_{\sigma \in S_k} \int_\gamma \omega_{\sigma(1)} \omega_{\sigma(2)} \cdots \omega_{\sigma(k)} = \int_\gamma \omega_1 \int_\gamma \omega_2 \cdots \int_\gamma \omega_k. \]
}

\item (Iterated Concatenation Formula)\label{SS:pairing-iter-coproduct}
\emph{Let $\gamma_1,\ldots,\gamma_r$ be paths in $\Gamma$ with $\gamma_{i+1}(0)=\gamma_i(1)$, and let $\omega_1,\ldots,\omega_k$ be polynomial $1$-forms. Then
      \[
      \int_{\gamma_1 \cdots \gamma_r} \omega_1 \cdots \omega_k =
        \sum_{\substack{g \in \Delta(k,r) \\ g = (n_1,\ldots,n_r)}}
          \Big(\int_{\gamma_1} \omega_1 \cdots \omega_{n_1}\Big)
          \Big(\int_{\gamma_2} \omega_{n_1+1} \cdots \omega_{n_1+n_2}\Big)
            \cdots
          \Big(\int_{\gamma_r} \omega_{n_1 + \cdots + n_{r-1} + 1} \cdots \omega_k\Big).
    \]
}

\end{itemize}
\end{theorem}

The concatenation formula together with the observation $\tint_{\beta\beta^{-1}} \omega_1\omega_2=0$ immediately yields the following \emph{conjugation formula} for iterated
integrals of length $2$.

\begin{theorem} Let $\alpha$ be a closed path in $\Gamma$, $\beta$ be a path with $\beta(1)=\alpha(0)$, and $\omega_1,\omega_2$ be polynomial $1$-forms. Then
\begin{align}
    \int_{\beta\alpha\beta^{-1}} \omega_1 \omega_2
      & = \int_\alpha \omega_1\omega_2
        + \left(\int_\beta\omega_1\int_\alpha\omega_2
          - \int_\alpha\omega_1\int_\beta\omega_2\right). \label{SS:pairing-conjugation}
\end{align}
\end{theorem}

For positive integers $k$ and $r$, denote by
\[ \Delta^+(k,r) \coloneqq \Set{f \in \Delta(k,r) | \text{$f$ is surjective}}. \]
By applying the iterated concatenation formula together with inclusion-exclusion, we see that combinatorial iterated integrals have the following useful nilpotence property:
\begin{theorem} \label{SS:pairing-nilpotence}
Let $\gamma_1,\ldots,\gamma_r$ be loops based at a point $v$, and let $\omega_1,\ldots,\omega_k$ be polynomial $1$-forms.
  Then
    \[
      \int_{(\gamma_1-1) \cdots (\gamma_r-1)} \omega_1 \cdots \omega_k =
        \sum_{\substack{g \in \Delta(k,r)^+ \\ g = (n_1,\ldots,n_r)}}
          \Big(\int_{\gamma_1} \omega_1 \cdots \omega_{n_1}\Big)
          \Big(\int_{\gamma_2} \omega_{n_1+1} \cdots \omega_{n_1+n_2}\Big)
            \cdots
          \Big(\int_{\gamma_r} \omega_{n_1 + \cdots + n_{r-1} + 1} \cdots \omega_k\Big).
    \]
  In particular, if $\omega_1,\ldots,\omega_k$ are tropical $1$-forms, interpreted as elements of $\rmH_1(\Gamma,\RR)$,
  \begin{equation}\label{eq:pairing-nilpotence.2}
    \int_{(\alpha_1 - 1) \cdots (\alpha_r - 1)} \omega_1 \cdots \omega_k =
      \begin{dcases}
        0 & r > k, \\
        \pairing{\omega_1,\alpha_1} \cdots \pairing{\omega_k,\alpha_k} & r = k.
      \end{dcases}
  \end{equation}
\end{theorem}

\begin{corollary}\label{SS:pairing-aug-vanish}
Let $J \coloneqq \ker(\ZZ\pi_1(\Gamma,v) \to \ZZ)$ be the augmentation ideal.
  For any $\alpha \in J^{\ell + 1}$, any $k \leq \ell$, and any
  $\omega_1,\ldots,\omega_k \in \rmH_1(\Gamma,\ZZ)$,
  \[ \int_\alpha \omega_1 \cdots \omega_k = 0. \]
  Thus, combinatorial iterated integration descends to a map
  \[ \int \colon \ZZ\pi_1(\Gamma,v)/J^{\ell + 1} \times T_\ell \rmH_1(\Gamma,\ZZ) \to \RR \]
  where $T_\ell \rmH_1(\Gamma,\ZZ) \coloneqq \bigoplus_{k = 0}^\ell \rmH_1(\Gamma,\ZZ)^{\otimes k}$ is the
  $\ell$\ssth~truncation of the tensor algebra.
\end{corollary}

\begin{proof}
  Any $\alpha \in J^{\ell + 1}$ can be written in the form
  \[
    \alpha =
      \sum_{r > \ell} \sum_{i_1,\ldots,i_r} c_{i_1 \ldots i_r}(\gamma_{i_1} - 1) \cdots (\gamma_{i_r} - 1)
  \]
  where all but finitely many of the $c_{i_1 \ldots i_r} \in \ZZ$ are $0$.
  It follows from~\eqref{eq:pairing-nilpotence.2} that $\int_\alpha$ is
  identically zero on $\rmH^1(\Gamma,\ZZ)^{\otimes k}$ for all $k \leq \ell$.
\end{proof}

The augmentation map $\epsilon\colon\ZZ\pi_1(\Gamma,v)\to\ZZ$ has kernel $J$ giving a descending filtration
\[\ZZ\pi_1(\Gamma,v)\supset J\supset J^2\supset\ldots.\]
Because $J/J^2\cong \rmsup{\pi_1(\Gamma,v)}{ab}\cong \rmH_1(\Gamma)$, the associated graded algebra can be identified as
\[ \gr_J(\ZZ\pi_1(\Gamma,v)) \cong \bigoplus\nolimits_{i = 0}^\infty (\rmH_1(\Gamma))^{\otimes i} \eqqcolon T(\rmH_1(\Gamma)). \]
Combinatorial iterated integrals can be interpreted as a bilinear map
\[ \int \colon \ZZ\pi_1(\Gamma,v) \times \gr_J(\ZZ\pi_1(\Gamma,v))   \rightarrow \RR.\]

\subsection{Integration algebras} \label{SS:pairing-algebra}
An \emph{integration algebra} over a ring $R$ is the data of $\big(A,\tint\big)$ where $A$ is an algebra over $\ZZ$ with nilradical $J$ inducing a filtration $A\supset J\supset J^2\supset\dots$ together with a bilinear map
\[
  \int \colon A \times \gr_J(A)\otimes\RR   \rightarrow \RR.
\]
Our natural example will be $\big(\ZZ\pi_1(\Gamma,v)/J^{\ell+1},\tint\big)$ where $\ell$ is a positive integer and $\tint$ is combinatorial iterated integration.

A \emph{morphism} of integration algebras over $R$,
\[ \varphi \colon \big(A, \tint\big) \to \big(A', \tint'\big) \]
is a $\ZZ$-algebra morphism
$\varphi \colon A \to A'$
preserving integration
in that
\[ \int_{\varphi(\gamma)} \gr_{J'}(\varphi)(\omega) = \int_\gamma \omega \]
for all $\gamma\in A$ and $\omega \in \gr_J(A)\otimes\RR$, where
$\gr_J(\varphi) \colon \gr_J(A)\otimes\RR \to \gr_{J'}(A')\otimes\RR$ is the
induced morphism on the associated graded algebras. Observe that because the filtration is induced by the nilradical, $\varphi$ must preserve the filtration.

Combinatorial iterated integration is nondegenerate in the following sense:
\begin{theorem}[Duality Theorem]\label{SS:pairing-duality}
For any pointed graph $(\Gamma,v)$ and each  $\ell \geq 0$, the maps
  \begin{equation}\label{eq:pairing-duality.1}
   \begin{aligned}
      T_\ell (\rmH_1(\Gamma)) \otimes_\ZZ \RR & \to \Hom_\ZZ(\ZZ \pi_1(\Gamma,v)/J^{\ell+1},\RR), &
      \RR \pi_1(\Gamma,v)/J^{\ell + 1} & \to \Hom_\ZZ(T_\ell(\rmH_1(\Gamma)),\RR), \\
        \omega_1\cdots\omega_k & \mapsto \int \omega_1\cdots\omega_k, &
        \gamma & \mapsto \int_\gamma,
    \end{aligned}
  \end{equation}
  are isomorphisms of vector spaces.
\end{theorem}

\begin{proof}
   The vector space $\RR \pi_1(\Gamma,v)/J^{\ell + 1}$ has a descending filtration
\[\RR \pi_1(\Gamma,v)/J^{\ell + 1}\supset J/J^{\ell + 1} \supset\dots\supset J^\ell/J^{\ell + 1}\]
while $T_\ell (\rmH_1(\Gamma)) \otimes_\ZZ \RR$ has an ascending filtration
\[(T_0 (\rmH_1(\Gamma)) \otimes_\ZZ \RR)\subset (T_1 (\rmH_1(\Gamma)) \otimes_\ZZ \RR) \subset\dots\subset (T_\ell (\rmH_1(\Gamma)) \otimes_\ZZ \RR).\]
  Consider the restriction of the combinatorial iterated integration
  \[ \int \colon \RR \pi_1(\Gamma,v)/J^{\ell + 1} \times T_\ell \rmH_1(\Gamma)  \otimes_\ZZ \RR \rightarrow\RR\]
to $J^{\ell}/J^{\ell+1}\times (\rmH_1(\Gamma) \otimes_\ZZ \RR)^{\otimes r}$. By the Nilpotence Property~\eqref{eq:pairing-nilpotence.2},
the pairing is $0$ if $r<\ell$ and is nondegenerate if $r=\ell$. Indeed, if $r=\ell$, the pairing factors as
\[J^r/J^{r+1}\times (\rmH_1(\Gamma)\otimes\RR)^{\otimes r}\cong (\rmH_1(\Gamma)\otimes\RR)^{\otimes r}\times (\rmH_1(\Gamma)\otimes\RR)^{\otimes r}\]
where it is the $r$-fold tensor product of the usual nondegenerate cycle pairing. From this upper-triangular structure of the combinatorial iterated integrals,
nondegeneracy follows.
\end{proof}

\section{Automorphism groups of integration algebras} \label{S:automorphisms}

In order to state our pointed Torelli theorem, we will need to study the automorphism groups of integration algebras. Let us consider $\rmH_1(\Gamma,\ZZ)$ as an abelian group equipped with the cycle pairing.

\begin{defn} The \emph{cyclic automorphism group} $\Aut_{\cyc}\big(\rmH_1(\Gamma,\ZZ),\pairing{\wc,\wc}\big)$ is the group of isometries
\[\phi\colon (\rmH_1(\Gamma,\ZZ),\pairing{\wc,\wc})\to (\rmH_1(\Gamma,\ZZ),\pairing{\wc,\wc}))\]
that induce the identity on $\rmH_1(\Gamma,\FF_2)$. The \emph{cyclic automorphism group} $\Aut_{\cyc}\big(\ZZ\pi_1(\Gamma,v)/J^3,\tint\big)$ is the group of isomorphisms of integration algebras
\[\phi\colon \big(\ZZ\pi_1(\Gamma,v)/J^3,\tint\big)\to \big(\ZZ\pi_1(\Gamma,v)/J^3,\tint\big)\]
that induce the identity on $\rmH_1(\Gamma,\FF_2)$.
\end{defn}

The group $\Aut_{\cyc}\big(\rmH_1(\Gamma,\ZZ),\pairing{\wc,\wc}\big)$ is always nontrivial because of the presence of multiplication by $-1$. In some sense, it is this automorphism that prevents the Torelli theorem for graphs \cite{CV10} from determining a $2$-connected graph up to isomorphism instead of $2$-isomorphism. Indeed, it foils an inductive approach to the Torelli theorem by preventing one from gluing isomorphisms.

We can give a complete description of  $\Aut_{\cyc}\big(\rmH_1(\Gamma,\ZZ),\pairing{\wc,\wc}\big)$. A connected graph can be written as an iterated $1$-point union of its $2$-connected components. That is, there are $2$-connected subgraphs $\Gamma_1,\dots,\Gamma_m$ such that there is a sequence of graphs $\Delta_1,\dots,\Delta_m=\Gamma$ with $\Delta_1=\Gamma_1$ and $\Delta_{i+1}=\Delta_i\vee_{u_{i+1},v_{i+1}} \Gamma_{i+1}$ where $\vee$ denotes a one-point union formed by identifying $u_{i+1}\in |\Delta_i|$ and $v_{i+1}\in|\Gamma_{i+1}|$. In this case, there is an orthogonal direct sum decomposition $\rmH_1(\Gamma,\ZZ)\cong \oplus_i \rmH_1(\Gamma_i,\ZZ)$.

Before we describe the cyclic automorphism groups, we recall some facts from \cite{CV10} about isometries $\phi\colon\rmH_1(\Gamma,\ZZ)\to \rmH_1(\Gamma',\ZZ)$ where $\Gamma$ and $\Gamma'$ are bridgeless graphs. A set of edges $S\subset E(\Gamma)$ is said to be a $C1$-set of $\Gamma$ if $\Gamma(S)$, the contraction of edges away from $S$, is a cycle and $\Gamma\setminus S$ has no bridges. The set of all $C1$-sets is denoted by $\operatorname{Set}^1(\Gamma)$. The $C1$-sets partition the edges of the graph, and every cyclic subgraph of $\Gamma$ can be written as a union of $C1$-sets. By \cite[Lemma~3.3.1]{CV10}, any $C1$-set arises as the intersection of two cyclic subgraphs in $\Gamma$. Any element of $\rmH_1(\Gamma,\ZZ)\subset \rmC_1(\Gamma,\ZZ)$ can be decomposed as an integer sum of $C1$-sets. For each $S\in \operatorname{Set}^1(\Gamma)$, there is a choice of orientation on the edges in $S$ such that when we write the sum of these oriented edges $e_S=\sum_{e\in S} e\in \rmC_1(\Gamma,\ZZ)$, for every cycle $C\in \rmH_1(\Gamma,\ZZ)$, we are able to write
\[C=\sum_{S\in \operatorname{Set}^1(\Gamma)} r_S(C)e_S\]
for some $r_S(C)\in\ZZ$. It is proven in \cite[Section~3.3]{CV10} that an isometry $\phi\colon\rmH_1(\Gamma,\ZZ)\to\rmH_1(\Gamma',\ZZ)$ induces a bijection $\beta\colon \operatorname{Set}^1(\Gamma)\to \operatorname{Set}^1(\Gamma')$ such that for any $C\in \rmH_1(\Gamma,\ZZ)$
\[r_{\beta(S)}(\phi(C))=\pm r_S(C).\]

\begin{prop} \label{p:homologyautomorphism} Let $\Gamma$ be a connected loopless graph. Let $\Gamma_1,\dots,\Gamma_m$ be the $2$-connected components of $\Gamma$. Then
\[\Aut_{\cyc}\big(\rmH_1,\Gamma,\ZZ,\pairing{\wc,\wc}\big)\cong \left(\ZZ/2\ZZ\right)^m\]
where $\sigma\in \left(\ZZ/2\ZZ\right)^m$ acts on the summand $\rmH_1(\Gamma_i,\ZZ)$ as multiplication by $(-1)^{\sigma_i}$.
\end{prop}

\begin{proof}
Let $\phi\colon \rmH_1(\Gamma,\ZZ)\to \rmH_1(\Gamma,\ZZ)$ be an cyclic automorphism. Produce the map
$\beta\colon\operatorname{Set}^1(\Gamma)\to \operatorname{Set}^1(\Gamma)$ as above. Since $\phi$ induces the identity on  $\rmH_1(\Gamma,\FF_2)$, for any cycle $C$, $\beta$ permutes the $C1$-sets supporting $C$, and $r_{\beta(S)}(\phi(C))=\pm r_S(C)$. By writing any $C1$-set $S$ as the intersection of two cyclic subgraphs, we see that $\beta$ is the identity map.
Write
\[r_S(\phi(C))=(-1)^{\sigma_{S,C}} r_S(C)\]
for $\sigma_{S,C}\in\ZZ/2\ZZ$. Now, given a primitive simple cycle $C\in\rmH_1(\Gamma,\ZZ)$, for $\phi(C)$ to be a cycle, we must have $\sigma_{S_1,C}=\sigma_{S_2,C}$ for all $C1$-sets $S_1,S_2$ supporting $C$. Therefore, for every primitive simple cycle $C$, $\phi(C)=\pm C$.

Because  $\Aut_{\cyc}\big(\rmH_1(\Gamma,\ZZ),\pairing{\wc,\wc}\big)$ decomposes as a direct product over $2$-connected components, we need only verify that if $\Gamma$ is $2$-connected, then $\Aut_{\cyc}\big(\rmH_1(\Gamma,\ZZ),\pairing{\wc,\wc}\big)\cong \{\pm 1\}$. Let $\phi$ be a cyclic automorphism of $\Gamma$. Because $\phi$ takes a primitive simple cycle $C$ to $\pm C$, and $\rmH_1(\Gamma,\ZZ)$ has a basis of primitive simple cycles, we can decompose
\[\rmH_1(\Gamma,\ZZ)\cong \rmH_1(\Gamma,\ZZ)^+\oplus \rmH_1(\Gamma,\ZZ)^-\]
where $\phi$ acts on the summands by $+1$ and $-1$, respectively. By the arguments above, the support of the cycles in $\rmH_1(\Gamma,\ZZ)^+$ is disjoint from the support of the cycles in $\rmH_1(\Gamma,\ZZ)^-$.
The partition of $E(\Gamma)$ into the edges supporting cycles in $\rmH_1(\Gamma,\ZZ)^+$ and $\rmH_1(\Gamma,\ZZ)^-$
 contradicts the connectedness of the matroid of $\Gamma$ \cite[Proposition~4.1.7]{Oxl11}.
\end{proof}

\begin{prop} \label{p:cpaautomorphism}
The natural map  $\Aut_{\cyc}\big(\ZZ\pi_1(\Gamma,v)/J^3,\tint\big)\to \Aut_{\cyc}\big(\rmH_1(\Gamma,\ZZ),\pairing{\wc,\wc}\big)$ is injective.
\end{prop}

\begin{proof}
Let $\phi\in \Aut_{\cyc}\big(\ZZ\pi_1(\Gamma,v)/J^3,\tint\big)$.
Pick free group generators $\gamma_1,\dots,\gamma_g$ for $\pi_1(\Gamma,v)$ such that their underlying cycle classes $[\gamma_i]\in\rmH_1(\Gamma,\ZZ)$ are each contained in a $2$-connected component of $\Gamma$. We claim that the automorphism $\phi$ is determined by $\phi([\gamma_1]),\dots,\phi([\gamma_g])$.
Write $\phi([\gamma_i])=(-1)^{\sigma_i}[\gamma_i]$.
Now, we know
\[\phi(\gamma_i)=(-1)^{\sigma_i}\gamma_i+\sum_{j,k} a_{ijk}(\gamma_j-1)(\gamma_k-1)\]
for integers $a_{ijk}$.
Let $\omega_1,\dots,\omega_g$ be a basis of $\rmH_1(\Gamma,\RR)$ dual to $[\gamma_1],\dots,[\gamma_g]$. Hence,
$\phi(\omega_i)=(-1)^{\sigma_i}\omega_i$.
Then,
\begin{align*}
\int_{\gamma_i} \omega_j\omega_k&=\int_{\phi(\gamma_i)} \phi(\omega_j)\phi(\omega_k)\\
&=(-1)^{\sigma_j+\sigma_k}\left((-1)^{\sigma_i}\int_{\gamma_i}\omega_j\omega_k+a_{ijk}\right).
\end{align*}
Consequently,
\[a_{ijk}=\left((-1)^{\sigma_j+\sigma_k}-(-1)^{\sigma_i}\right)\int_{\gamma_i}\omega_j\omega_k.\]
\end{proof}

The group $\Aut_{\cyc}\big(\ZZ\pi_1(\Gamma,v)/J^3,\tint\big)$ is somewhat mysterious. As we will see in Section~\ref{S:hyperelliptic},
if $\Gamma$ is a $2$-connected hyperelliptic graph and $v$ is a fixed point of a hyperelliptic involution, then $\Aut_{\cyc}\big(\ZZ\pi_1(\Gamma,v)/J^3,\tint\big)$ is nontrivial and hence, by Proposition~\ref{p:cpaautomorphism},  isomorphic to $\ZZ/2\ZZ$. Conjecture~\ref{conj:hyperelliptic} is the converse of this statement.

We expect that a positive resolution of this conjecture implies Conjecture \ref{conj:unipotenttorelli} and thus that $\big(\ZZ\pi_1(\Gamma,v)/J^3,\tint\big)$ and the harmonic volume of $\Gamma$ (defined below) are complete invariants of $(\Gamma,v)$.

\section{Integration Algebras and Base-points}\label{S:unipotent}

The integration algebra determines the base-point of a bridgeless graph $\Gamma$ with $g(\Gamma)\geq 2$ up to well-understood finite ambiguity. Let us first recall the Abel--Jacobi map of graphs following \cite{MZ-Jacobian,Baker-Norine}.

Let $\Gamma$ be a metric graph with base-point $v$. The cycle pairing gives an injection $\rmH_1(\Gamma,\ZZ)\to \rmH^1(\Gamma,\RR)\cong \rmH_1(\Gamma,\RR)^\vee$ by, for $C\in \rmH_1(\Gamma,\ZZ)$,
\[C\mapsto \left[\omega\mapsto \int_C \omega\right].\]
The Jacobian of $\Gamma$ is the cokernel of this map,
\[\sJ_1(\Gamma)\coloneqq \rmH^1(\Gamma,\RR)/\rmH_1(\Gamma,\ZZ).\]
Now, we recall the \emph{Abel--Jacobi map with respect to $v\in |\Gamma|$}, $\iota_v\colon |\Gamma|\to \sJ_1(\Gamma)$. We define
\[\iota_v(p)\coloneqq \left[\omega\mapsto \int_\alpha \omega\right].\]
where $\alpha$ is any path from $v$ to $p$.
This function is independent of the choice of $\alpha$.

The following is well-known \cite{Baker-Norine,BX11}:
\begin{theorem}The map $\iota_v$ is injective if and only if $\Gamma$ is bridgeless.
\end{theorem}

The Abel--Jacobi map is particularly straightforward when we work with dual bases. Let $C_1,\dots,C_g$ be a basis of $\rmH_1(\Gamma,\ZZ)$, and let $\omega_1,\dots,\omega_g$ be the dual basis of $\rmH_1(\Gamma,\RR)$ under $\pairing{\wc,\wc}$. Then, by using the inner product on $\rmH_1(\Gamma,\RR)$, we may view the Abel--Jacobi map as
\begin{align*}
\iota_v \colon |\Gamma|&\to \RR^g/\ZZ^g\\
p&\mapsto \left(\int_\alpha \omega_1,\dots,\int_\alpha \omega_g\right)
\end{align*}

\begin{defn} A {\em cycle-respecting isomorphism} from $\big(\ZZ\pi_1(\Gamma,v)/J^3,\tint\big)\to \big(\ZZ\pi_1(\Gamma,v')/J^3,\tint\big)$ is an isomorphism of integration algebras inducing the identity on $\rmH_1(\Gamma,\FF_2)$.
\end{defn}

\begin{theorem}\label{SS:basepoint} Let $\Gamma$ be a bridgeless metric graph with $g(\Gamma)\geq 2$. Let $v$ be a  point of $|\Gamma|$. Then the cycle-respecting isomorphism type of the integration algebra $\big(\ZZ\pi_1(\Gamma,v)/J^3,\tint\big)$ determines $v$ up to at most
$\left|\Aut_{\cyc}\big(\ZZ\pi_1(\Gamma,v)/J^3,\tint\big)\right|$ choices.
\end{theorem}

\begin{proof}
Let $v'$ be a point of $\Gamma$ such that there is a cycle-respecting isomorphism between $\big(\ZZ\pi_1(\Gamma,v)/J^3,\tint\big)$ and $\big(\ZZ\pi_1(\Gamma,v')/J^3,\tint\big)$. Then there are $\big|\Aut_{\cyc}\big(\ZZ\pi_1(\Gamma,v)/J^3,\tint\big)\big|$ such isomorphisms. We will show that the data of an isomorphism
$\phi\colon \big(\ZZ\pi_1(\Gamma,v)/J^3,\tint\big)\to \big(\ZZ\pi_1(\Gamma,v')/J^3,\tint\big)$
determines $\iota_{v'}(v)$. It will also follow from our arguments that if $\Aut_{\cyc}\big(\ZZ\pi_1(\Gamma,v)/J^3,\tint\big)$ is trivial, then $\iota_{v'}(v)=0$ and thus $v=v'$.

Let $\alpha$ be a path in $\Gamma$ from $v'$ to $v$. Let $C_1,\dots,C_g\in H_1(\Gamma,\ZZ)$ be cycles, each supported in a $2$-connected component of $\Gamma$, giving a basis. Let $\omega_1,\dots,\omega_g\in H_1(\Gamma,\RR)$ be the dual basis under $\pairing{\wc,\wc}$. Let $\gamma_1,\dots,\gamma_g$ be a free group basis for $\pi_1(\Gamma,v)$ such that $[\gamma_i]=C_i$. Write
$\phi([C_i])=(-1)^{\sigma_i}[C_i]$ for $\sigma_i\in \ZZ/2\ZZ$.

Now
\[\phi(\gamma_i)=(-1)^{\sigma_i}\alpha\gamma_i\alpha^{-1}+\sum_{j,k} a_{ijk} (\gamma_j-1)(\gamma_k-1)\]
for integers $a_{ijk}$. Pick $j\neq i$.
We have
\begin{align*}
\int_{\gamma_i} \omega_j\omega_i&=\int_{\phi(\gamma_i)} \phi(\omega_j)\phi(\omega_i)\\
&= (-1)^{\sigma_j}\int_{\alpha\gamma_i\alpha^{-1}} \omega_j\omega_i+(-1)^{\sigma_i+\sigma_j}a_{iji}
\end{align*}
By (\ref{SS:pairing-conjugation}),
\begin{align*}
\int_{\alpha\gamma_i\alpha^{-1}} \omega_j\omega_i&=\int_{\gamma_i} \omega_j\omega_i+\left(\int_{\alpha} \omega_j\int_{\gamma_i} \omega_i-\int_{\gamma_i} \omega_j\int_{\alpha}\omega_i\right)\\
&=\int_{\gamma_i} \omega_j\omega_i+\int_{\alpha} \omega_j
\end{align*}
Therefore,
\[\int_\alpha \omega_j=\left((-1)^{\sigma_j}-1\right)\int_{\gamma_i}\omega_j\omega_i-(-1)^{\sigma_i}a_{iji},
\]
Because $a_{iji}\in\ZZ$, the image of the Abel--Jacobi map is determined by the combinatorial iterated integrals.

In the case that $\sigma_i=0$ for all $i$, we have that $\int_\alpha \omega_j$ is always an integer. Hence, $\iota_{v'}(v)=0$.
\end{proof}

More can be said if $\Gamma$ is a $3$-connected graph. The Torelli theorem for metric graphs due to Caporaso--Viviani \cite[Theorem~4.1.9]{CV10} shows that the isometry class of $\rmH_1(\Gamma,\ZZ)$ determines $\Gamma$ up to the relation called tropical equivalence. Here, if $v$ is a $2$-valent vertex adjacent to edges $e_1$ and $e_2$, \emph{smoothing} $v$ replaces $e_1$, $e_2$, and $v$ by an edge
$e$ with $\ell(e)=\ell(e_1)+\ell(e_2)$. The operation inverse to smoothing is \emph{refinement}. Two bridgeless graphs $\Gamma,\Gamma'$ are \emph{tropically equivalent} if they have isometric refinements.

\begin{corollary} Let $(\Gamma,v)$ be a pointed metric graphs such that $\Gamma$ is $3$-connected. Then the integration algebra $\big(\ZZ\pi_1(\Gamma,v)/J^3,\tint\big)$ determines $\Gamma$ up to tropical equivalence and determines $v$ up to two possibilities.
\end{corollary}

\begin{proof}
Suppose that $(\Gamma,v)$ and $(\Gamma',v')$ are pointed graphs such that there is an isomorphism $\phi\colon
\big(\ZZ\pi_1(\Gamma,v)/J^3,\tint\big)\to \big(\ZZ\pi_1(\Gamma',v')/J^3,\tint\big)$. Then $\phi$ induces an isometry between $\rmH_1(\Gamma,\ZZ)$ and $\rmH_1(\Gamma',\ZZ)$. Because $\Gamma$ is $3$-connected, $\Gamma$ and $\Gamma'$ are tropically equivalent by the Torelli theorem for metric graphs. Moreover, the isometry $\phi$ is induced by this tropical equivalence, a fact that is implicit in \cite{CV10}. One sees it by noting that the proof of the Torelli theorem produces a bijection between $C1$-sets $\operatorname{Set}^1 \Gamma\to \operatorname{Set}^1 \Gamma'$. For $3$-connected graphs, there is a natural bijection $E(\Gamma)\cong \operatorname{Set}^1\Gamma$. Consequently, there is a unique 
cyclic bijection $\Phi\colon E(\Gamma)\to E(\Gamma')$. This cyclic bijection induces the tropical equivalence by the strong form of Whitney's $2$-isomorphism theorem for $3$-connected graphs \cite[Lemma~5.3.2]{Oxl11}.
Therefore, $\phi$ is a cycle-preserving isomorphism.  Since $\Gamma$ is $2$-connected, by Proposition~\ref{p:homologyautomorphism} and Proposition~\ref{p:cpaautomorphism}, $\left|\Aut_{\cyc}\big(\ZZ\pi_1(\Gamma,v)/J^3,\tint\big)\right|$ has at most two elements.
\end{proof}

\section{The Harmonic Volume of Pointed Graphs}\label{S:comp}

In this section, we introduce the \emph{harmonic volume}, an invariant of pointed graphs within a $2$-isomorphism class,
analogous to Bruno Harris's harmonic volume in algebraic geometry \cite{HarVol}.
It encodes the data of the integration algebra $\big(\ZZ\pi_1(\Gamma,v)/J^3,\tint\big)$ as an extension of
$\rmH_1(\Gamma,\ZZ)$. It is valued in a real torus, $\sJ_2(\Gamma)$ analogous
to an intermediate Jacobian.

Let $\Gamma$ be a connected, metric, bridgeless graph. Define
\begin{align*}
  \sJ_2(\Gamma)
    \coloneqq \frac{\Hom_\ZZ\big(\rmH_1(\Gamma,\ZZ),\rmH^1(\Gamma,\RR) \otimes \rmH^1(\Gamma,\RR)\big)}{%
        \Hom_\ZZ\big(\rmH_1(\Gamma,\ZZ),\rmH_1(\Gamma,\ZZ) \otimes \rmH_1(\Gamma,\ZZ)\big)}.
\end{align*}
Here, the quotient comes from viewing
$\rmH_1(\Gamma,\ZZ) \otimes \rmH_1(\Gamma,\ZZ)$ as a lattice in
$\rmH^1(\Gamma,\RR) \otimes \rmH^1(\Gamma,\RR)$ via the inclusion
$\rmH_1(\Gamma,\ZZ) \hookrightarrow \rmH^1(\Gamma,\RR)\cong \rmH_1(\Gamma,\RR)^\vee$ induced by the cycle pairing
$\pairing{\wc,\wc}$.

A \emph{pointed rigged graph} in the $2$-isomorphism class of $\Gamma$ is a triple $(\Gamma',v',\phi')$ where $(\Gamma',v')$ is a pointed metric, bridgeless graph, and $\phi'\colon\rmH_1(\Gamma,\ZZ)\to \rmH_1(\Gamma',\ZZ)$ is an isometry. Let $W(\Gamma)$ be the set of pointed rigged graphs $2$-isomorphic to $\Gamma$.

We now produce the harmonic volume of $(\Gamma',v')$.
Consider the following exact sequence where $J$ is the augmentation ideal of $\ZZ\pi_1(\Gamma',v')/J^3$:
\[
  \begin{tikzcd}
    0 \rar & J^2/J^3 \rar & J/J^3 \rar & J/J^2  \rar & 0.
  \end{tikzcd}
\]
By using the isomorphisms
\[J/J^2\cong \rmH_1(\Gamma',\ZZ),\quad J^2/J^3\cong \rmH_1(\Gamma',\ZZ)\otimes \rmH_1(\Gamma',\ZZ),\]
we can write the diagram as
\[  \begin{tikzcd}
    0 \rar & \rmH_1(\Gamma',\ZZ)\otimes \rmH_1(\Gamma',\ZZ) \rar & J/J^3 \rar["\operatorname{[\cdot ]}"] & \rmH_1(\Gamma',\ZZ) \arrow[bend right,swap]{l}{\sigma}  \rar & 0 \\
 \end{tikzcd}
\]
 where $[\cdot]$ maps $\gamma-1\in J/J^3$ to its underlying cycle $[\gamma]$.
Given any section $\sigma$ as above, we define  the harmonic volume as the map
 $\mu\colon \rmH_1(\Gamma,\ZZ)\to (\rmH_1(\Gamma,\RR)\otimes\rmH_1(\Gamma,\RR))^\vee$ given by
 \[C\mapsto \left[\omega_1\otimes\omega_2\mapsto \int_{\sigma(\phi(C))} \phi(\omega_1)\phi(\omega_2)\right].\]
 The class of $\mu$ in $\sJ_2(\Gamma)$, which we denote by $\nu(\Gamma,v)$,  is independent of the choice of lift $\sigma$. Indeed,
 any two sections $\sigma_1,\sigma_2$ differ by a map
\[\sigma_1-\sigma_2\colon \rmH_1(\Gamma',\ZZ)\to \rmH_1(\Gamma',\ZZ)\otimes \rmH_1(\Gamma',\ZZ).\]

\begin{theorem}\label{SS:comp-ext-determine} Let $(\Gamma_1,v_1,\phi_1),(\Gamma_2,v_2,\phi_2)\in W(\Gamma)$. There is equality of harmonic volumes, $\nu_\Gamma(\Gamma_1,v_1,\phi_1)=\nu_\Gamma(\Gamma_2,v_2,\phi_2)$ if and only if  $\phi_2\circ\phi_1^{-1}\colon \rmH_1(\Gamma_1,\ZZ)\to \rmH_1(\Gamma_2,\ZZ)$ lifts to an isomorphism
\[\phi\colon \big(\ZZ\pi_1(\Gamma_1,v_1)/J_1^3,\tint\big)\to \big(\ZZ\pi_1(\Gamma_2,v_2)/J_2^3,\tint\big).\]
\end{theorem}

\begin{proof}
  Pick lifts $\sigma_i\colon \rmH_1(\Gamma_i,\ZZ)\to J_i/J_i^3$ to define the harmonic volumes.   Because $\nu_\Gamma(\Gamma_1,v_1,\phi_1)=\nu_\Gamma(\Gamma_2,v_2,\phi_2)$, we may modify $\sigma_2$ such that for any $C\in \rmH_1(\Gamma,\ZZ)$ and $\omega_1,\omega_2\in
  \rmH_1(\Gamma,\RR)$,
  \[\int_{\sigma_1(\phi_1(C))}\phi_1(C)\phi_1(C)=\int_{\sigma_2(\phi_2(C))} \phi_2(C)\phi_2(C).\]
  We use $\sigma_i$ and $\phi_i$ to write isomorphisms
  \begin{align*}
  \ZZ\pi_1(\Gamma_i,v_i)/J_1^3
  &\cong \ZZ\oplus \rmH_1(\Gamma_i,\ZZ)\oplus
  \rmH_1(\Gamma_i,\ZZ)\otimes \rmH_1(\Gamma_i,\ZZ)\\
  &\cong \ZZ\oplus \rmH_1(\Gamma,\ZZ)\oplus  \rmH_1(\Gamma,\ZZ)\otimes \rmH_1(\Gamma,\ZZ).
  \end{align*}
  Define $\phi$ to be the composition of this isomorphism for $i=1$ and its inverse for $i=2$.
  \[\ZZ\pi_1(\Gamma_1,v_1)/J_1^3\to \ZZ\oplus \rmH_1(\Gamma,\ZZ)\oplus
  \rmH_1(\Gamma,\ZZ)\otimes \rmH_1(\Gamma,\ZZ) \to \ZZ\pi_1(\Gamma_2,v_2)/J_2^3.\]
  In other words, this is the unique isomorphism that restricts to $\phi$ when taking the quotient by $J_i^2$ and which intertwines $\sigma_1$ and $\sigma_2$. It is straightforward to verify that this is an isomorphism of integration algebras.
\end{proof}

\begin{corollary} Let $\Gamma$ be a $2$-connected graph. Let $(\Gamma',v',\phi')\in W(\Gamma)$. Then $\Aut_{\cyc}\big(\ZZ\pi_1(\Gamma',v')/J^3,\tint\big)$ is non-trivial if and only if $\nu_\Gamma(\Gamma',v',\phi')$ is $2$-torsion.
\end{corollary}

\begin{proof}
The non-triviality of $\Aut_{\cyc}\big(\ZZ\pi_1(\Gamma',v')/J^3,\tint\big)$ is equivalent to multiplication by $-1$ on $\rmH_1(\Gamma',\ZZ)$ lifting to some non-trivial automorphism of $\big(\ZZ\pi_1(\Gamma',v')/J^3,\tint\big)$.
This occurs if and only if $\nu_\Gamma(\Gamma',v',\phi')=\nu_\Gamma(\Gamma',v',-\phi')$. Since $\nu_\Gamma(\Gamma',v',-\phi')=-\nu_\Gamma(\Gamma',v',\phi')$, the conclusion follows.
\end{proof}

Conjecture~\ref{conj:hyperelliptic} is equivalent to $\nu(\Gamma',v',\phi)$ being $2$-torsion if and only $(\Gamma',v')$ is hyperelliptic.

\section{Hyperelliptic graphs}\label{S:hyperelliptic}

\subsection{Background on hyperelliptic graphs}

Hyperelliptic graphs  were studied by Baker and Norine in \cite{BNHyperelliptic}. For $2$-edge connected graphs $\Gamma$, hyperellipticity is equivalent to the existence of a {\em hyperelliptic involution} $\iota\colon\Gamma\to\Gamma$ such that $\Gamma/\iota$ is a tree. If points $u,w\in |\Gamma|$ are interchanged by $\iota$, they are said to be hyperelliptically conjugate. If $(\Gamma,v)$ is a pointed graph such that there exists a hyperelliptic involution fixing $v$, then we say the pair $(\Gamma,v)$ is hyperelliptic.

Let $(\Gamma,v)$ be a $2$-connected hyperelliptic pointed graph. Then, $\Aut_{\cyc}\big(\ZZ\pi_1(\Gamma,v)/J^3,\tint\big)$ is non-trivial with non-trivial element induced by $\iota$. Indeed, $\iota$ acts as the identity on $\rmH_1(\Gamma,\FF_2)$ as can be seen as follows. Let $T$ be the tree $\Gamma/\iota$, $\pi\colon\Gamma\to T$ be the quotient map,
and $\partial T$ be the image of the fixed points of $T$. Cyclic subgraphs of $\Gamma$ are the preimages under
$\pi$ of simple paths between points in $\partial T$. These are fixed by $\iota$. This can also be seen directly by an explicit description of  $\rmH_1(\Gamma,\ZZ)$ as in \cite[Section~4]{KK20}. Moreover, $\iota$ acts as multiplication by $-1$ on $\rmH_1(\Gamma,\ZZ)$.

\subsection{Potential-theoretic Criterion for Hyperellipticity}

We present Theorem~\ref{t:hyperellipticity}, a criterion for hyperellipticity which was part of an unsuccessful attempt to prove Conjecture~\ref{conj:hyperelliptic} but which we think is of  interest. It would be worthwhile to compare this criterion to Corey's criterion for graphs of hyperelliptic type \cite{Corey19}.

\begin{defn} A {\em piecewise linear function} on $\Gamma$ is a function $\psi\colon |\Gamma|\to \RR$ that is linear on every edge of some refinement of $\Gamma$.
Attached to $\psi$ is a degree $0$ polynomial $1$-form $d\psi$ whose value on an oriented edge $e$ is $p_e(t) dt$. Given a degree $0$ polynomial $1$-form $\omega$, we may define $d^*\omega$ to be the formal real combination of points of $|\Gamma|$ where
\[d^*\omega=\sum_{q\in|\Gamma|}\left(\sum_e p_e(t)\right)q\]
where the inner sum is over the edges $e$ directed away from $q$.
The {\em Laplacian} is defined by $\Delta(\psi)\coloneqq d^*d\psi$.
\end{defn}

It is easily seen that a degree $0$ polynomial $1$-form $\omega$ is equal to $d\psi$ for some piecewise linear $\psi$ exactly when $\int_{\gamma} \omega=0$ for all closed paths $\gamma$.

\begin{defn} Let $(\Gamma,v)$ be a pointed graph.
Given a piecewise linear function $\psi$ and $C_1,C_2\in \rmH_1(\Gamma,\RR)$, we define the pairing with respect to $\psi$ by
\[\pairing{C_1,C_2}_\psi=\int_\gamma \psi \omega_{C_2}\]
where $\gamma$ is a loop based at $v$ with $[\gamma]=C_1$. This is independent of the choice of $\gamma$.
\end{defn}

The above pairing is symmetric. Of particular interest will be the case where $\psi$ is a piecewise linear function with $\Delta(\psi)=u-w$ for points $u,w\in |\Gamma|$. In this case, $\psi$ is defined up to a real constant and $\pairing{\wc,\wc}_\psi$ is well-defined up to taking the sum with a real multiple of the usual cycle pairing $\pairing{\wc,\wc}$. This pairing naturally occurs in performing certain iterated integrals.

The following is a consequence of integration by parts.
\begin{lemma}  \label{l:psiintegral}
Let $\psi$ be a piecewise linear function on $\Gamma$.
Then,
\[\pairing{C_1,C_2}_\psi=\int_\gamma (d\psi)\omega_{C_2}+\left(\psi(\gamma(0))\pairing{C_1,C_2}\right)\]
where $\gamma\colon[0,1]\to|\Gamma|$ is any loop with $\gamma(0)=v$ and $[\gamma]=C_1$.
\end{lemma}

Before we begin the proof of the hyperellipticity criterion, recall that a banana graph is a graph with two vertices, no loops, and at least two edges between the vertices.

\begin{theorem} \label{t:hyperellipticity} Let $\Gamma$ be a $2$-connected metric graph. Let $u,w\in |\Gamma|$. Let $\psi$ be a piecewise linear function on $\Gamma$ such that $\Delta(\psi)=u-w$. Then \(\Gamma\) is hyperelliptic making $u$ and $w$ conjugate if and only if $\pairing{\wc,\wc}_\psi$ is a scalar multiple of $\pairing{\wc,\wc}$.
\end{theorem}

\begin{proof}
First suppose that there is an involution $\iota\colon\Gamma\to\Gamma$ with $\iota(u)=w$ such that $\Gamma/\iota$ is a tree. By adding a constant to $\psi$, we may suppose $\psi(v)=-\psi(u)$. Then $\psi(u)<0<\psi(w)$. Because $\Delta(\iota^*\psi)=w-u=\Delta(-\psi)$
and $\iota^*\psi(u)=-\psi(u)$, $\iota^*\psi=-\psi$. Since for any $\omega\in \rmH_1(\Gamma,\RR)$, $\iota^*\omega=-\omega$,
\begin{align*}
\pairing{\omega_1,\omega_2}_\psi&=\pairing{\iota^*\omega_1,\iota^*\omega_2}_{\iota^*\psi}\\
&=\pairing{-\omega_1,-\omega_2}_{-\psi}\\
&=-\pairing{\omega_1,\omega_2}_\psi,
\end{align*}
it follows that $\pairing{\omega_1,\omega_2}_\psi=0$.

Now suppose that $\pairing{\wc,\wc}_\psi$ is a scalar multiple of $\pairing{\wc,\wc}$. By adding a constant to $\psi$, we may suppose $\pairing{\wc,\wc}_\psi=0$. We will write $\Gamma$ as the union of two trees. Let
\[\Gamma_-\coloneqq \overline{\psi^{-1}((-\infty,0))},\quad \Gamma_+\coloneqq \overline{\psi^{-1}((0,\infty))},
\quad \Gamma_0=\psi^{-1}(0).\]
We define the boundary of $\Gamma_{\pm}$ to be $\partial\Gamma_{\pm}=\Gamma_{\pm}\cap \Gamma_0$.

We first claim that $\Gamma_-$ is a tree. First observe that $\Gamma_-$ contains no cycles because if $C$ were a simple cycle in $\Gamma_-$, then $\pairing{\omega_C,\omega_C}_\psi<0$. Now, suppose by way of contradiction that $\Gamma_-$ is not connected. Then it has a component $K$ not containing $u$. If $K^\circ=K\setminus \Gamma_0$ denotes the interior of $K$, $\Delta(\psi|_{K^\circ})=0$ and $\psi$ cannot obtain its minimum on $K^\circ$. Since $\psi|_K$ obtains its maximum on $\partial K$, it must be a constant. Because $K$ is a component of $\Gamma_-$, $\psi$ must be negative at some point of $K$, and thus $K$ must be a component of $\Gamma$. This contradiction gives the claim. Similarly, $\Gamma_+$ is a tree. Observe that $\Delta(\psi|_{\Gamma_-})$ is strictly negative on $\partial \Gamma_-$ and
$\Delta(\psi|_{\Gamma_+})$ is strictly positive on $\partial \Gamma_+$.

We claim that $\Gamma_0$ is a finite set of points. Suppose to the contrary that $K$ is a component of $\Gamma_0$
containing an edge. Now, $K$ intersects in $\Gamma_-$ and $\Gamma_+$ in at most one point. Indeed, if
$|K\cap \Gamma_-|\geq 2$, $\Gamma_-\cup K$ contains a simple cycle $C$ for which $\pairing{\omega_C,\omega_C}_\psi<0$ since one can form a cycle from paths $\gamma_1,\gamma_2$ from a point in $\Gamma_-^{\circ}$ to points in $K$ together with a path in $K$. Because $\Gamma$ is $2$-connected, $K$ must intersect both $\Gamma_-$ and $\Gamma_+$. Let $w_-=\Gamma_-\cap K$ and $w_+=\Gamma_+\cap K$. If $w_-=w_+$, then $w_-$ would be a disconnecting vertex of $\Gamma$, hence $w_-\neq w_+$.
 Because $\psi$ is constant on $K$,
$\Delta(\psi|_K)(w_-)=0$.
Then
\[\Delta(\psi)(w_-)=\Delta(\psi|_{\Gamma_-})(w_-)+\sum_{K}\Delta(\psi|_{K})(w_-)<0\]
where the sum is over the components $K$ of $\Gamma_0$ containing $w_-$.
This contradiction shows that no such $K$ exists.

Consequently $\Gamma=\Gamma_-\cup\Gamma_+$. Moreover, each point of $\Gamma=\Gamma_-\cap\Gamma_+$ is $2$-valent in $\Gamma$. Otherwise, $\Gamma_-$ or $\Gamma_+$ would fail to be a tree.

On a $\psi$-increasing path in $\Gamma_-$, the slopes of $\psi$ are non-increasing. Indeed,
given any point $z$ in $\Gamma_-$, we can form a path from $z$ by following edges along which $\psi$ is non-increasing. This path must terminate at the unique minimum $u$ of $\psi$. Since $\Gamma_-$ is a tree, there is only one path from $z$ to $u$. Consequently, at any $z\in |\Gamma_-|$, there is at most one edge along which $\psi$ is non-increasing. Therefore, because $\Delta(\psi)(z)=0$ for $z\neq u$, the sum of the positive slopes of $\psi$ along edges incident to $z$ is equal to the opposite of the unique negative slope of an edge incident to $z$. Consequently, the slope of $\psi$ on a $\psi$-increasing edge from $z$ is at most the opposite of the slope of $\psi$ on the $\psi$-non-increasing edge from $z$. From this the claim about slopes follows. Similarly, on a $\psi$-increasing path in $\Gamma_+$, the slopes of $\psi$ are non-decreasing.

We claim that there is an involution $\iota$ of $\Gamma$ interchanging $(\Gamma_-,u)$ and $(\Gamma_+,w)$ and fixing $\Gamma_0$. We induct on the number of vertices of $\Gamma$. If $\Gamma$ has two vertices $u,w$, $\Gamma$ is a banana graph, and we're done. For the inductive step, put a graph structure on $\Gamma$ such that $\psi$ is linear on each edge. We may suppose that $\Gamma$ has no $2$-valent vertices except possibly $u$ and $w$. Because $\Gamma_-$ and $\Gamma_+$ are both trees, there are no vertices $z$ for which $\psi(z)=0$. Now, let $z_-$ be a vertex for which $|\psi(z_-)|$ is minimal. By replacing $\psi$ with $-\psi$, we may suppose that $\psi(z_-)<0$. Let $e_1,\dots,e_k$ be the edges adjacent to $z_-$ on which $\psi$ increases from $z_-$. Given any pair of edges $e_i,e_j$, we may follow them to get paths $\gamma_i$ and $\gamma_j$ from $z_-$ that meet in $\Gamma_+$ at a point $z_{ij}$. Observe that if $z_i$ and $z_j$ are the first vertices of $\gamma_i$ and $\gamma_j$ in $\Gamma_+$, then $e_i=z_-z_i$, $e_j=z_-z_j$ and
\[\psi(z_{i})\geq -\psi(z_-),\ \psi(z_j)\geq -\psi(z_-).\]
Let $C=[\gamma_i]-[\gamma_j]$.
Observe that the slopes of $\psi$ along $\gamma_i$ (resp.~$\gamma_j$) are non-decreasing, and they are constant from $z_-$ to $z_i$ (resp. $z_j)$ since there
are no other vertices to encounter.
Therefore,
\[\int_{e_i} \psi\omega_C\geq 0,\quad \int_{e_j} \psi\omega_C\geq 0.\]
The integral of $\psi\omega_C$ over any other edge of $C$ must be positive. Because $\int_C \psi\omega_C=0$, we must have $C=e_i-e_j$ and thus $z_{ij}=z_i=z_j$
Therefore, all the $\psi$-increasing edges from $z_-$ meet at the same point $z_+=z_{ij}$,
yielding a banana subgraph. We replace that subgraph by a single edge $e'=z_-z_+$ to obtain a graph $\Gamma'$. We
define a piecewise linear function $\psi'$ on $\Gamma'$. It is defined to equal $\psi$ away from $e'$ and to be
linear on $e'$. The length of $e'$ is chosen such that the slope of $\psi'$ on $e'$ is equal to the sum of the slopes of $\psi$
on $e_1,\dots,e_k$.
Consequently, $\Delta(\psi')(z_-)=\Delta(\psi')(z_+)=0$. By smoothing  $z_-$ and $z_+$, we see $\Gamma'$ has two
fewer vertices than $\Gamma$. Moreover, $\pairing{\wc,\wc}_{\psi'}$ vanishes. Indeed, let $C'$ is a cycle in $\Gamma'$. If it is
supported away from $e'$, then
$\pairing{\omega_{C'},\omega_{C'}}_{\psi'}=\pairing{\omega_{C'},\omega_{C'}}_{\psi}=0$. If it
contains $e'$, we replace $e'$ by $e_1$ to get a cycle $C$ in $\Gamma$. Then
$\pairing{\omega_{C'},\omega_{C'}}_{\psi'}=\pairing{\omega_{C},\omega_{C}}_{\psi}=0$.
By induction, there is an involution $\iota'\colon\Gamma'\to\Gamma'$ interchanging $u$ and $v$ whose quotient is a tree. By replacing the involution on $e'$ by the corresponding involution of the banana graph, we obtain the desired involution $\iota\colon\Gamma\to\Gamma$.
\end{proof}

\bibliography{main}
\bibliographystyle{amsalpha}
\end{document}